\documentclass[10pt]{article}

\usepackage{amssymb}
\usepackage{eucal}
\usepackage{amsmath}
\usepackage{amsthm}
\usepackage{amscd}
\usepackage{graphics}
\usepackage{graphicx}
\usepackage{amsfonts}
\usepackage{latexsym}
\usepackage[all]{xy}

\numberwithin{equation}{section}
\newtheorem{thm}{\bf Theorem}[section]

\theoremstyle{remark}

\title{ON THE STABILITY PROBLEM FOR THE $\mathfrak{so}(5)$ FREE RIGID BODY}

\author{IOAN CA\c SU\\
{\small Department of Mathematics, West University of Timi\c soara,}\\
{\small Bd. V. P\^{a}rvan, No. 4, 300223 Timi\c soara, Romania}\\ 
{\small E-mail: casu@math.uvt.ro}\\
}

\date{}

\begin{document}

\maketitle




\begin{abstract}
In the general case of the $\mathfrak{so}(n)$ free rigid body we will give a list of integrals of motion, which generate the set of Mishchenko's integrals.
In the case of
$\mathfrak{so}(5)$ we prove that there are fifteen coordinate type Cartan subalgebras which on a regular adjoint orbit give fifteen Weyl group orbits of equilibria. 
These coordinate type
Cartan subalgebras are the analogues of the three axes of equilibria
for the classical rigid body on $\mathfrak{so}(3)$.
The nonlinear stability and instability of these equilibria is analyzed.
In addition to these equilibria there are ten other continuous families of equilibria. 
\end{abstract}


\section{Introduction}

The purpose of this paper is to study the geometry and the dynamics of the free rigid body on the Lie algebra $\mathfrak{so}(5)$. The general free rigid body on $\mathfrak{so}(n)$ has been studied as a completely integrable system in the classical works \cite{Mishchenko70}, \cite{Manakov76}, \cite{Ratiu80}. The general case and more extensively the even case $n=2m$ have been analyzed by Feh\'{e}r and Marshall \cite{FeMa03}, Spiegler 
\cite{Spiegler04}, where a method for determining a certain class of equilibria is given. They also study the stability of these equilibria using the Energy-Casimir method. As mentioned in \cite{Spiegler04}, the odd case $n=2m+1$ is significantly different and only some indications on stability results have been given.

In \cite{noi} it has been proved that this class of equilibria comes as the intersection of coordinate type Cartan subalgebras with regular adjoint orbits. Consequently these equilibria are unions of Weyl orbits and they correspond to long-short axis type of equilibria known from the dynamics on $\mathfrak{so}(3)$. 

The stability of these equilibria, for the case of $\mathfrak{so}(4)$, has been studied in \cite{noi}, using Williamson normal form \cite{BoFo04}. 
For the case of $\mathfrak{so}(4)$ there have been discovered two new continuous families of equilibria on every regular adjoint orbit. Also, in 
\cite{noi} it has been shown that these are nonlinearly stable as a family, that is, if a solution of the $\mathfrak{so}(4)$-free rigid body equation starts near an equilibrium on such a curve, at any later time it will stay close to this curve but in the direction of the curve itself it may drift.

First of all, in this paper we will give a list of integrals of motion for the general case of the $\mathfrak{so}(n)$ free rigid body, which proves to be a set of generating functions for the Mishchenko's quadratic integrals of motion and which will play a crucial role in the approach of the stability problem using energy methods. 

For the case of $\mathfrak{so}(5)$ rigid body we found fifteen families of coordinate type Cartan subalgebras, which intersected with a regular adjoint orbit give 120=15$\times$8 equilibria corresponding to long-short axis type of equilibria (8 being the cardinal of the Weyl group of $\mathfrak{so}(5)$). We also found ten continuous families of equilibria on every regular adjoint orbit. As in the case of $\mathfrak{so}(4)$, these two types of equilibria give all the equilibria of $\mathfrak{so}(5)$ free rigid body. 

The stability problem is studied using the linearization method and Arnold's method, which is equivalent to Energy-Casimir method \cite{BiPu07}. In order to do this, we explicitly compute the integrals of motion using Mishchenko's and Manakov's methods and we point out the need for the additional constants of motion introduced in Section 2 of this paper. Extended indications on how Arnold's method applies for each studied equilibria are given.

\section{A new family of integrals of motion for the free rigid body on $\mathfrak{so}(n)$}

The equations of the rigid body on $\mathfrak{so}(n)$ are given by
\begin{equation}\label{rigideq}
\dot M=[M,\Omega], \end{equation} where $\Omega\in \mathfrak{so}(n)$,
$M=\Omega J+J\Omega\in \mathfrak{so}(n)$ with $J=\operatorname{diag}(\lambda_i)$, a real constant diagonal matrix satisfying $ \lambda _i + \lambda _j \geq 0 $, for all $i, j=1, \ldots, n$, $i \neq j $ (see, for example, \cite{Ratiu80}).
Note that $M=[m_{ij}]$ and $\Omega=[ \omega_{ij}]$ determine each other if and only if $\lambda_i + \lambda _j >0$ since $m_{ij} = ( \lambda_i+ \lambda_j) \omega_{ij}$ which physically means that the rigid body is not concentrated on a lower dimensional subspace of $\mathbb{R}^n$. 

It is well known and easy to verify that equations \eqref{rigideq} are Hamiltonian relative to the minus Lie-Poisson bracket 
\begin{equation}
\label{LP_n}
\{F,G\}(M) : = \frac{1}{2}\operatorname{Trace}(M [\nabla F(M), \nabla G(M)]),
\end{equation}
and the Hamiltonian function
\begin{equation}
\label{ham_n}
H(M) := - \frac{1}{4} \operatorname{Trace}(M \Omega).
\end{equation}
Here $F,G, H \in C ^{\infty}(\mathfrak{so}(n))$ and the gradient is taken relative to the Ad-invariant inner product
\begin{equation}
\label{ip_n}
\left\langle X, Y \right\rangle : = - \frac{1}{2}\operatorname{Trace}(XY), 
\qquad X, Y \in \mathfrak{so}(n)
\end{equation}
which identifies $(\mathfrak{so}(n))^*$ with $\mathfrak{so}(n)$. This means that $\dot{F} = \{F, H\}$ for all $F \in C^{\infty}(\mathfrak{so}(n))$, where $\{\cdot , \cdot \}$ is given by \eqref{LP_n} and $H $ by \eqref{ham_n}, if and only if \eqref{rigideq} holds.
Note that the linear isomorphism $X \in \mathfrak{so}(n) \mapsto XJ + JX 
\in \mathfrak{so}(n)$ is self-adjoint relative to the inner product \eqref{ip_n} and thus $\nabla H(M) =  \Omega$.

In what follows, we assume that $\lambda_i$ are all distinct. As stated in \cite{Mishchenko70}, the system \eqref{rigideq} admits a sequence $m_r$ ($r=1,2,3,\dots$) of integrals of motion, all of them depending quadratically on the angular momentum $M$. It is easy to compute an explicit form for these integrals \cite{MoPi} as follows:
$$m_r(M)=\sum\limits_{\stackrel{i,k=1}{i<k}}^n\frac{\lambda_i^r-\lambda_k^r}{\lambda_i^2-\lambda_k^2}m_{ik}^2.$$
Next we will find an additional family of $n$ integrals of motion for the system \eqref{rigideq}, which have simple and elegant expressions and which prove to generate Mishchenko's integrals of motion. The complete result is contained in the following
\begin{thm}
\label{generators}
The functions
$$F_i(M)=\sum\limits_{\stackrel{k=1}{k\not= i}}^{n}\frac{1}{\lambda_i^2-\lambda_k^2}m_{ik}^2,~~~i=1,2,\dots,n$$
are integrals of motion for the system \eqref{rigideq}; moreover, the following identities hold for each $r$:
$$m_r(M)=\sum\limits_{i=1}^n\lambda_i^rF_i(M).$$
\end{thm}

\noindent {\bf Proof.} It is easy to see that
\begin{align*}
m_r(M)&=\sum\limits_{\stackrel{i,k=1}{i<k}}^n\frac{\lambda_i^r}{\lambda_i^2-\lambda_k^2}m_{ik}^2-\sum\limits_{\stackrel{i,k=1}{i<k}}^n\frac{\lambda_k^r}{\lambda_i^2-\lambda_k^2}m_{ik}^2\\
&=\sum\limits_{\stackrel{i,k=1}{i<k}}^n\frac{\lambda_i^r}{\lambda_i^2-\lambda_k^2}m_{ik}^2+\sum\limits_{\stackrel{i,k=1}{k<i}}^n\frac{\lambda_i^r}{\lambda_i^2-\lambda_k^2}m_{ki}^2\\
&=\sum_{i=1}^n\lambda_i^r\left(\sum\limits_{\stackrel{k=1}{k>i}}^n\frac{1}{\lambda_i^2-\lambda_k^2}m_{ik}^2+\sum\limits_{\stackrel{k=1}{k<i}}^n\frac{1}{\lambda_i^2-\lambda_k^2}m_{ik}^2\right)\\
&=\sum_{i=1}^n\lambda_i^r\cdot\sum\limits_{\stackrel{k=1}{k\not= i}}^n\frac{1}{\lambda_i^2-\lambda_k^2}m_{ik}^2\\
&=\sum_{i=1}^n\lambda_i^r F_i(M).
\end{align*}
Now, writing the above identities for $r=1,2,\dots,n$ we obtain a Vandermonde-type linear system, which solved leads to the conclusion that the functions $F_i$ are linear combinations of Mishchenko's integrals of motion $m_r$ ($r=1,2,\dots,n$) and thus are integrals of motion for our system \eqref{rigideq}. \rule{0.5em}{0.5em}

\section{The free rigid body on the Lie algebra $\mathfrak{so}(5)$}

We begin by making some considerations on $\mathfrak{so}(5)$, which is the Lie algebra of the subgroup 
$\operatorname{SO}(5)=\{A\in \mathfrak{gl}(5, \mathbb{R}) 
\mid A^tA=I_5,\det (A)=1\}$ of the special linear Lie group
$\operatorname{SL}(5,\mathbb{R})$.

We choose as basis of $\mathfrak{so}(5)$ the matrices
$$E_1=\left[\begin{array}{ccccc}
0&0&0&0&0\\
0&0&-1&0&0\\
0&1&0&0&0\\
0&0&0&0&0\\
0&0&0&0&0
\end{array}\right];~E_2=\left[\begin{array}{ccccc}
0&0&1&0&0\\
0&0&0&0&0\\
-1&0&0&0&0\\
0&0&0&0&0\\
0&0&0&0&0\end{array}\right];$$
$$E_3=\left[\begin{array}{ccccc}
0&-1&0&0&0\\
1&0&0&0&0\\
0&0&0&0&0\\
0&0&0&0&0\\
0&0&0&0&0\end{array}\right];
E_4=\left[\begin{array}{ccccc}
0&0&0&1&0\\
0&0&0&0&0\\
0&0&0&0&0\\
-1&0&0&0&0\\
0&0&0&0&0\end{array}\right];$$
$$E_5=\left[\begin{array}{ccccc}
0&0&0&0&0\\
0&0&0&1&0\\
0&0&0&0&0\\
0&-1&0&0&0\\
0&0&0&0&0\end{array}\right];E_6=\left[\begin{array}{ccccc}
0&0&0&0&0\\
0&0&0&0&0\\
0&0&0&1&0\\
0&0&-1&0&0\\
0&0&0&0&0\end{array}\right];$$
$$
E_7=\left[\begin{array}{ccccc}
0&0&0&0&1\\
0&0&0&0&0\\
0&0&0&0&0\\
0&0&0&0&0\\
-1&0&0&0&0\end{array}\right];E_8=\left[\begin{array}{ccccc}
0&0&0&0&0\\
0&0&0&0&1\\
0&0&0&0&0\\
0&0&0&0&0\\
0&-1&0&0&0\end{array}\right];$$
$$E_9=\left[\begin{array}{ccccc}
0&0&0&0&0\\
0&0&0&0&0\\
0&0&0&0&1\\
0&0&0&0&0\\
0&0&-1&0&0\end{array}\right];
E_{10}=\left[\begin{array}{ccccc}
0&0&0&0&0\\
0&0&0&0&0\\
0&0&0&0&0\\
0&0&0&0&1\\
0&0&0&-1&0\end{array}\right]$$
and hence we represent $\mathfrak{so}(5)$ as
\begin{equation}
\label{so4_representation}
\mathfrak{so}(5)=\left\{\, M=\left|\, \left[\begin{array}{ccccc}
0&-x_3&x_2&y_1&z_1\\
x_3&0&-x_1&y_2&z_2\\
-x_2&x_1&0&y_3&z_3\\
-y_1&-y_2&-y_3&0&z_4\\
-z_1&-z_2&-z_3&-z_4&0
\end{array}\right] \, \right|\, x_1,x_2,x_3,y_1,y_2,y_3,z_1,z_2,z_3,z_4\in \mathbb{R}\right\}.
\end{equation}
Relative to this basis, the Lie algebra structure of $\mathfrak{so}(5)$ is given by the following table
\medskip
\begin{center}
\begin{tabular}{|c|c|c|c|c|c|c|c|c|c|c|}
\hline
$[\cdot,\cdot]$&$E_1$&$E_2$&$E_3$&$E_4$&$E_5$&$E_6$&$E_7$&$E_8$&$E_9$&$E_{10}$\\
\hline
$E_1$&0&$E_3$&$-E_2$&0&$E_6$&$-E_5$&0&$E_9$&$-E_8$&0\\
\hline
$E_2$&$-E_3$&0&$E_1$&$-E_6$&0&$E_4$&$-E_9$&0&$E_7$&0\\
\hline
$E_3$&$E_2$&$-E_1$&0&$E_5$&$-E_4$&0&$E_8$&$-E_7$&0&0\\
\hline
$E_4$&0&$E_6$&$-E_5$&0&$E_3$&$-E_2$&$-E_{10}$&0&0&$E_7$\\
\hline
$E_5$&$-E_6$&0&$E_4$&$-E_3$&0&$E_1$&0&$-E_{10}$&0&$E_8$\\
\hline
$E_6$&$E_5$&$-E_4$&0&$E_2$&$-E_1$&0&0&0&$-E_{10}$&$E_9$\\
\hline
$E_7$&0&$E_9$&$-E_8$&$E_{10}$&0&0&0&$E_3$&$-E_2$&$-E_4$\\
\hline
$E_8$&$-E_9$&0&$E_7$&0&$E_{10}$&0&$-E_3$&0&$E_1$&$-E_5$\\
\hline
$E_9$&$E_8$&$-E_7$&0&0&0&$E_{10}$&$E_2$&$-E_1$&0&$-E_6$\\
\hline
$E_{10}$&0&0&0&$E_7$&$-E_8$&$-E_9$&$E_4$&$E_5$&$E_6$&0\\
\hline
\end{tabular}
\end{center}
\medskip
In the basis $\{E_1, \ldots, E_{10}\}$, the matrix of the Lie-Poisson
structure \eqref{LP_n} is 
\begin{equation}
\label{LP}
\Gamma_-=\left[\begin{array}{cccccccccc}
0&-x_3&x_2&0&-y_3&y_2&0&-z_3&z_2&0\\
x_3&0&-x_1&y_3&0&-y_1&z_3&0&-z_1&0\\
-x_2&x_1&0&-y_2&y_1&0&-z_2&z_1&0&0\\
0&-y_3&y_2&0&-x_3&x_2&z_4&0&0&-z_1\\
y_3&0&-y_1&x_3&0&-x_1&0&z_4&0&-z_2\\
-y_2&y_1&0&-x_2&x_1&0&0&0&z_4&-z_3\\
0&-z_3&z_2&-z_4&0&0&0&-x_3&x_2&y_1\\
z_3&0&-z_1&0&-z_4&0&x_3&0&-x_1&y_2\\
-z_2&z_1&0&0&0&-z_4&-x_2&x_1&0&y_3\\
0&0&0&z_1&z_2&z_3&-y_1&-y_2&-y_3&0
\end{array}\right].
\end{equation}
Since $\hbox{rank}\,\mathfrak{so}(5) = 2$, there are two functionally independent Casimir functions which are given respectively by
\[
C_1(M) := -\frac{1}{4}\operatorname{Trace}(M^2)
= \frac{1}{2}\left(\sum\limits\limits_{i=1}^3 x_i^2+\sum\limits\limits_{i=1}^3 y_i^2+\sum\limits\limits_{i=1}^4 z_i^2 \right)
\]
and
\[
C_2(M) := \frac{1}{8}\operatorname{Trace}(M^4)=
\]
$$=\frac{1}{8}\left[(x_2^2+x_3^2+y_1^2+z_1^2)^2+(x_1^2+x_3^2+y_2^2+z_2^2)^2+(x_1^2+x_2^2+y_3^2+z_3^2)^2+\right.$$
$$+(y_1^2+y_2^2+y_3^2+z_4^2)^2+(z_1^2+z_2^2+z_3^2+z_4^2)^2+2(y_1z_1+y_2z_2+y_3z_3)^2+$$
$$+2(x_2y_3-x_3y_2-z_1z_4)^2+2(x_3y_1-x_1y_3-z_2z_4)^2+2(x_1y_2-x_2y_1-z_3z_4)^2+$$
$$+2(x_3z_2-x_2z_3-y_1z_4)^2+2(x_1z_3-x_3z_1-y_2z_4)^2+2(x_2z_1-x_1z_2-y_3z_4)^2+$$
$$\left. +2(x_1x_2-y_1y_2-z_1z_2)^2+2(x_1x_3-y_1y_3-z_1z_3)^2+2(x_2x_3-y_2y_3-z_2z_3)^2 \right].$$
Thus the generic adjoint orbits are the level sets
$$
\operatorname{Orb}_{c_1c_2}(M)=(C_1\times C_2)^{-1}(c_1,c_2), 
\qquad (c_1,c_2) \in \mathbb{R}^2.
$$

In all that follows we will denote by $\operatorname{Orb}_{c_1;c_2}$ the regular adjoint orbit $\operatorname{Orb}_{c_1c_2}$,
where $c_1>0,c_2>0$ and $2c_2>c_1^2>c_2$.

Using the Lie bracket table in the chosen basis given above, it is immediately seen that the coordinate type Cartan subalgebras of $\mathfrak{so}(5)$ are
$\mathfrak{t}_k,~~1\leq k\leq 15$, where
$$\mathfrak{t}_1:=\operatorname{span}(E_3,E_6)=\left\{\left. M_{a,b}^1:=\left[\begin{array}{ccccc}
0&a&0&0&0\\
-a&0&0&0&0\\
0&0&0&b&0\\
0&0&-b&0&0\\
0&0&0&0&0\end{array}\right]\, \right| \,a,b\in \mathbb{R}
 \right\},$$
$$\mathfrak{t}_2:=\operatorname{span}(E_6,E_8)=\left\{\left. M_{a,b}^2:=\left[\begin{array}{ccccc}
0&0&0&0&0\\
0&0&0&0&a\\
0&0&0&-b&0\\
0&0&b&0&0\\
0&-a&0&0&0\end{array}\right]\, \right| \,a,b\in \mathbb{R}
 \right\},$$
$$\mathfrak{t}_3:=\operatorname{span}(E_6,E_7)=\left\{\left. M_{a,b}^3:=\left[\begin{array}{ccccc}
0&0&0&0&a\\
0&0&0&0&0\\
0&0&0&-b&0\\
0&0&b&0&0\\
-a&0&0&0&0\end{array}\right]\, \right| \,a,b\in \mathbb{R}
 \right\},$$
$$\mathfrak{t}_4:=\operatorname{span}(E_5,E_7)=\left\{\left. M_{a,b}^4:=\left[\begin{array}{ccccc}
0&0&0&0&a\\
0&0&0&b&0\\
0&0&0&0&0\\
0&-b&0&0&0\\
-a&0&0&0&0\end{array}\right]\, \right| \,a,b\in \mathbb{R}
 \right\},$$ 
$$\mathfrak{t}_5:=\operatorname{span}(E_1,E_7)=\left\{\left. M_{a,b}^:=\left[\begin{array}{ccccc}
0&0&0&0&a\\
0&0&-b&0&0\\
0&b&0&0&0\\
0&0&0&0&0\\
-a&0&0&0&0\end{array}\right]\, \right| \,a,b\in \mathbb{R}
 \right\},$$
$$\mathfrak{t}_6:=\operatorname{span}(E_2,E_5)=\left\{\left. M_{a,b}^6:=\left[\begin{array}{ccccc}
0&0&a&0&0\\
0&0&0&-b&0\\
-a&0&0&0&0\\
0&b&0&0&0\\
0&0&0&0&0\end{array}\right]\, \right| \,a,b\in \mathbb{R}
 \right\},$$
$$\mathfrak{t}_7:=\operatorname{span}(E_5,E_9)=\left\{\left. M_{a,b}^7:=\left[\begin{array}{ccccc}
0&0&0&0&0\\
0&0&0&a&0\\
0&0&0&0&-b\\
0&-a&0&0&0\\
0&0&b&0&0\end{array}\right]\, \right| \,a,b\in \mathbb{R}
 \right\},$$
$$\mathfrak{t}_8:=\operatorname{span}(E_1,E_{10})=\left\{\left. M_{a,b}^8:=\left[\begin{array}{ccccc}
0&0&0&0&0\\
0&0&a&0&0\\
0&-a&0&0&0\\
0&0&0&0&-b\\
0&0&0&b&0\end{array}\right]\, \right| \,a,b\in \mathbb{R}
 \right\},$$
$$\mathfrak{t}_9:=\operatorname{span}(E_1,E_4)=\left\{\left. M_{a,b}^9:=\left[\begin{array}{ccccc}
0&0&0&a&0\\
0&0&-b&0&0\\
0&b&0&0&0\\
-a&0&0&0&0\\
0&0&0&0&0\end{array}\right]\, \right| \,a,b\in \mathbb{R}
 \right\},$$
$$\mathfrak{t}_{10}:=\operatorname{span}(E_2,E_8)=\left\{\left. M_{a,b}^{10}:=\left[\begin{array}{ccccc}
0&0&a&0&0\\
0&0&0&0&-b\\
-a&0&0&0&0\\
0&0&0&0&0\\
0&b&0&0&0\end{array}\right]\, \right| \,a,b\in \mathbb{R}
 \right\},$$
$$\mathfrak{t}_{11}:=\operatorname{span}(E_4,E_8)=\left\{\left. M_{a,b}^{11}:=\left[\begin{array}{ccccc}
0&0&0&a&0\\
0&0&0&0&-b\\
0&0&0&0&0\\
-a&0&0&0&0\\
0&b&0&0&0\end{array}\right]\, \right| \,a,b\in \mathbb{R}
 \right\},$$
$$\mathfrak{t}_{12}:=\operatorname{span}(E_3,E_{10})=\left\{\left. M_{a,b}^{12}:=\left[\begin{array}{ccccc}
0&a&0&0&0\\
-a&0&0&0&0\\
0&0&0&0&0\\
0&0&0&0&b\\
0&0&0&-b&0\end{array}\right]\, \right| \,a,b\in \mathbb{R}
 \right\},$$
$$\mathfrak{t}_{13}:=\operatorname{span}(E_3,E_9)=\left\{\left. M_{a,b}^{13}:=\left[\begin{array}{ccccc}
0&a&0&0&0\\
-a&0&0&0&0\\
0&0&0&0&-b\\
0&0&0&0&0\\
0&0&b&0&0\end{array}\right]\, \right| \,a,b\in \mathbb{R}
 \right\},$$
$$\mathfrak{t}_{14}:=\operatorname{span}(E_4,E_9)=\left\{\left. M_{a,b}^{14}:=\left[\begin{array}{ccccc}
0&0&0&a&0\\
0&0&0&0&0\\
0&0&0&0&-b\\
-a&0&0&0&0\\
0&0&b&0&0\end{array}\right]\, \right| \,a,b\in \mathbb{R}
 \right\},$$
$$\mathfrak{t}_{15}:=\operatorname{span}(E_2,E_{10})=\left\{\left. M_{a,b}^{15}:=\left[\begin{array}{ccccc}
0&0&a&0&0\\
0&0&0&0&0\\
-a&0&0&0&0\\
0&0&0&0&b\\
0&0&0&-b&0\end{array}\right]\, \right| \,a,b\in \mathbb{R}
 \right\}.$$ 

The intersection of a regular adjoint orbit and a coordinate
type Cartan subalgebra has eight elements which represents a Weyl group orbit. Specifically, we have the following result:

\begin{thm} \label{Weil_group_orbit}
The following equalities hold:
$$\mathfrak{t}_k\cap \operatorname{Orb}_{c_1;c_2}=\left\{
M_{a,b}^k,M_{-a,-b}^k,M_{b,a}^k,M_{-b,-a}^k,M_{-a,b}^k,M_{a,-b}^k,M_{b,-a}^k,M_{-b,a}^k
\right\}$$
for all $k$, $1\leq k\leq 15$,
where
\begin{equation}
\label{values_a_b}
\left\{
\begin{aligned}
a&=\sqrt{c_1+\sqrt{2c_2-c_1^2}}\\
b&=\sqrt{c_1-\sqrt{2c_2-c_1^2}}.
\end{aligned}
\right.
\end{equation}
\end{thm}

\noindent {\bf Proof.} 
Let $M_{\alpha,\beta}^1\in \mathfrak{t}_1\cap \operatorname{Orb}_{c_1;c_2}$. 
Then $M_{\alpha,\beta}^1\in\mathfrak{t}_1$, $c_1=C_1(M_{\alpha,\beta}^1)=\frac{1}{2}\left(\alpha^2+\beta^2\right)$ and $c_2=C_2(M_{\alpha,\beta}^1)=\frac{1}{4}\left(\alpha^4+\beta^4\right)$.
We obtain the system
$$\left\{\begin{array}{l}
\alpha^2+\beta^2=2c_1\\
\alpha^4+\beta^4=4c_2\end{array}\right.$$
which leads immediately to the result. For $\mathfrak{t}_k\cap \operatorname{Orb}_{c_1;c_2}$, $2\leq k\leq 15$, we proceed in a similar manner.
\rule{0.5em}{0.5em}
\medskip


\section{Equilibria for the $\mathfrak{so}(5)$-rigid body}

We will work from now on with a generic $\mathfrak{so}(5)$-rigid body, that is, $\lambda_i+ \lambda_j>0$ for $i \neq j$ and all $\lambda_i $ are distinct. 
The relationship between $\Omega=[\omega_{ij}]\in \mathfrak{so}(5)$ and the matrix $M\in \mathfrak{so}(5)$ in the representation \eqref{so4_representation} is hence given by
$$(\lambda_3+\lambda_2)\omega_{32}=x_1; \quad 
(\lambda_1+\lambda_3)\omega_{13}=x_2; \quad
(\lambda_2+\lambda_1)\omega_{21}=x_3;$$
$$(\lambda_1+\lambda_4)\omega_{14}=y_1; \quad 
(\lambda_2+\lambda_4)\omega_{24}=y_2; \quad 
(\lambda_3+\lambda_4)\omega_{34}=y_3;$$
$$(\lambda_1+\lambda_5)\omega_{15}=z_1; \quad 
(\lambda_2+\lambda_5)\omega_{25}=z_2; \quad 
(\lambda_3+\lambda_5)\omega_{35}=z_3; \quad 
(\lambda_4+\lambda_5)\omega_{45}=z_4.$$
and thus the equations of motion \eqref{rigideq} are equivalent to the system
\begin{equation}\label{soM}
\left\{\begin{array}{l}
\vspace{.1cm}
\dot x_1=(\lambda_2-\lambda_3)\left[\frac{y_2y_3}{(\lambda_2+\lambda_4)(\lambda_3+\lambda_4)}-\frac{x_2x_3}{(\lambda_1+\lambda_2)(\lambda_1+\lambda_3)}+\frac{z_2z_3}{(\lambda_2+\lambda_5)(\lambda_3+\lambda_5)}\right]\\\\
\dot x_2=(\lambda_3-\lambda_1)\left[\frac{y_1y_3}{(\lambda_1+\lambda_4)(\lambda_3+\lambda_4)}-\frac{x_1x_3}{(\lambda_1+\lambda_2)(\lambda_2+\lambda_3)}+\frac{z_1z_3}{(\lambda_1+\lambda_5)(\lambda_3+\lambda_5)}\right]\\\\
\dot x_3=(\lambda_1-\lambda_2)\left[\frac{y_1y_2}{(\lambda_1+\lambda_4)(\lambda_2+\lambda_4)}-\frac{x_1x_2}{(\lambda_2+\lambda_3)(\lambda_1+\lambda_3)}+\frac{z_1z_2}{(\lambda_1+\lambda_5)(\lambda_2+\lambda_5)}\right]\\\\
\dot y_1=(\lambda_1-\lambda_4)\left[\frac{x_2y_3}{(\lambda_1+\lambda_3)(\lambda_3+\lambda_4)}-\frac{x_3y_2}{(\lambda_1+\lambda_2)(\lambda_2+\lambda_4)}-\frac{z_1z_4}{(\lambda_1+\lambda_5)(\lambda_4+\lambda_5)}\right]\\\\
\dot y_2=(\lambda_2-\lambda_4)\left[\frac{x_3y_1}{(\lambda_1+\lambda_2)(\lambda_1+\lambda_4)}-\frac{x_1y_3}{(\lambda_2+\lambda_3)(\lambda_3+\lambda_4)}-\frac{z_2z_4}{(\lambda_2+\lambda_5)(\lambda_4+\lambda_5)}\right]\\\\
\dot y_3=(\lambda_3-\lambda_4)\left[\frac{x_1y_2}{(\lambda_2+\lambda_3)(\lambda_2+\lambda_4)}-\frac{x_2y_1}{(\lambda_1+\lambda_3)(\lambda_1+\lambda_4)}-\frac{z_3z_4}{(\lambda_3+\lambda_5)(\lambda_4+\lambda_5)}\right]\\\\
\dot z_1=(\lambda_1-\lambda_5)\left[\frac{x_2z_3}{(\lambda_1+\lambda_3)(\lambda_3+\lambda_5)}-\frac{x_3z_2}{(\lambda_1+\lambda_2)(\lambda_2+\lambda_5)}+\frac{y_1z_4}{(\lambda_1+\lambda_4)(\lambda_4+\lambda_5)}\right]\\\\
\dot z_2=(\lambda_2-\lambda_5)\left[
\frac{x_3z_1}
{(\lambda_1+\lambda_2)(\lambda_1+\lambda_5)}-\frac{x_1z_3}{(\lambda_2+\lambda_3)(\lambda_3+\lambda_5)}+\frac{y_2z_4}{(\lambda_2+\lambda_4)(\lambda_4+\lambda_5)}\right]\\\\
\dot z_3=(\lambda_3-\lambda_5)\left[\frac{x_1z_2}{(\lambda_2+\lambda_3)(\lambda_2+\lambda_5)}-\frac{x_2z_1}{(\lambda_1+\lambda_3)(\lambda_1+\lambda_5)}+\frac{y_3z_4}{(\lambda_3+\lambda_4)(\lambda_4+\lambda_5)}\right]\\\\
\dot z_4=(\lambda_4-\lambda_5)\left[-\frac{y_1z_1}{(\lambda_1+\lambda_4)(\lambda_1+\lambda_5)}-\frac{y_2z_2}{(\lambda_2+\lambda_4)(\lambda_2+\lambda_5)}-\frac{y_3z_3}{(\lambda_3+\lambda_4)(\lambda_3+\lambda_5)}\right].
\end{array}\right.
\end{equation}
The Hamiltonian \eqref{ham_n} has in this case the expression
\begin{align*}
H(M)&=-\frac{1}{4}\hbox{Trace}(M\Omega)\\
&=\frac{1}{2}\left(\frac{1}{\lambda_2+\lambda_3}x_1^2+\frac{1}{\lambda_1+\lambda_3}x_2^2+\frac{1}{\lambda_1+\lambda_2}x_3^2+\frac{1}{\lambda_1+\lambda_4}y_1^2+\frac{1}{\lambda_2+\lambda_4}y_2^2+\frac{1}{\lambda_3+\lambda_4}y_3^2+\right.\\
&+\left.\frac{1}{\lambda_1+\lambda_5}z_1^2+\frac{1}{\lambda_2+\lambda_5}z_2^2+\frac{1}{\lambda_3+\lambda_5}z_3^2+\frac{1}{\lambda_4+\lambda_5}z_4^2\right).
\end{align*}
The Hamiltonian nature of system \eqref{soM} can be checked in this case directly, writing 
$$(\dot{x}_1, \dot{x}_2, \dot{x}_3, \dot{y}_1, \dot{y}_2, \dot{y}_3,\dot{z}_1,\dot{z}_2,\dot{z}_3,\dot{z}_4)^{\sf T} = \Gamma_- (\nabla H)^{\sf T},$$ where the Poisson structure $\Gamma_- $ is given by \eqref{LP}.

\begin{thm} \label{echilibre} 
If $\mathcal{E}$ denotes the set of the equilibrium points
of \eqref{soM}, then
${\cal E}=\left(\bigcup\limits_{k=1}^{15}\mathfrak{t}_k\right)\cup \left(\bigcup\limits_{l=1}^{10}\mathfrak{s}_l\right)$,
where $\mathfrak{s}_l,~1\leq l\leq 10$, are the three dimensional vector subspaces given by\\

\medskip

$\mbox{\fontsize{10}{10}\selectfont $
\mathfrak{s}_{1,2}:=\operatorname{span}_{\mathbb{R}}\left\{\left(\frac{1}{\lambda_1+\lambda_4}E_1 \pm \frac{1}{\lambda_2+\lambda_3}E_4\right),\,\left(\frac{1}{\lambda_2+\lambda_4}E_2\pm \frac{1}{\lambda_1+\lambda_3}E_5\right), \,\left(\frac{1}{\lambda_3+\lambda_4}E_3\pm \frac{1}{\lambda_1+\lambda_2}E_6\right) \right\};$}$

\medskip

$\mbox{\fontsize{10}{10}\selectfont $
\mathfrak{s}_{3,4}:=\operatorname{span}_{\mathbb{R}}\left\{\left(\frac{1}{\lambda_4+\lambda_5}E_1 \pm \frac{1}{\lambda_2+\lambda_3}E_{10}\right),\,\left(\frac{1}{\lambda_3+\lambda_5}E_5\pm \frac{1}{\lambda_2+\lambda_4}E_9\right), \,\left(\frac{1}{\lambda_2+\lambda_5}E_6\mp \frac{1}{\lambda_3+\lambda_4}E_8\right) \right\};$}$

\medskip

$\mbox{\fontsize{10}{10}\selectfont $
\mathfrak{s}_{5,6}:=\operatorname{span}_{\mathbb{R}}\left\{\left(\frac{1}{\lambda_4+\lambda_5}E_2 \pm \frac{1}{\lambda_1+\lambda_3}E_{10}\right),\,\left(\frac{1}{\lambda_1+\lambda_5}E_6\pm \frac{1}{\lambda_3+\lambda_4}E_7\right), \,\left(\frac{1}{\lambda_3+\lambda_5}E_4\mp \frac{1}{\lambda_1+\lambda_4}E_9\right) \right\};$}$

\medskip

$\mbox{\fontsize{10}{10}\selectfont $
\mathfrak{s}_{7,8}:=\operatorname{span}_{\mathbb{R}}\left\{\left(\frac{1}{\lambda_1+\lambda_5}E_1 \pm \frac{1}{\lambda_2+\lambda_3}E_7\right),\,\left(\frac{1}{\lambda_2+\lambda_5}E_2\pm \frac{1}{\lambda_1+\lambda_3}E_8\right), \,\left(\frac{1}{\lambda_3+\lambda_5}E_3\pm \frac{1}{\lambda_1+\lambda_2}E_9\right) \right\};$}$

\medskip

$\mbox{\fontsize{10}{10}\selectfont $
\mathfrak{s}_{9,10}:=\operatorname{span}_{\mathbb{R}}\left\{\left(\frac{1}{\lambda_4+\lambda_5}E_3 \pm \frac{1}{\lambda_1+\lambda_2}E_{10}\right),\,\left(\frac{1}{\lambda_2+\lambda_5}E_4\pm \frac{1}{\lambda_1+\lambda_4}E_8\right), \,\left(\frac{1}{\lambda_1+\lambda_5}E_5\mp \frac{1}{\lambda_2+\lambda_4}E_7\right) \right\}.$}$
\end{thm}

\noindent {\bf Proof.} The proof follows after a long, but straightforward computation. \rule{0.5em}{0.5em}

\section{Constants of motion and nonlinear stability}
In this section we study the nonlinear stability of the
equilibrium states ${\cal E}\cap \operatorname{Orb}_{c_1;c_2}$ for the dynamics \eqref{soM} on a generic adjoint orbit.

Using Mishchenko's method \cite{Mishchenko70}, \cite{Ratiu80}, \cite{MoPi}, we obtain the following additional independent constants of the motion for Eqs. \eqref{soM}:
\begin{align*}
K_1(M)&=-\frac{1}{4}\operatorname{Trace}\sum\limits_{p=0}^3 J^{p}MJ^{3-p}\Omega\\
&=\frac{1}{2}[(\lambda_2^2+\lambda_3^2)x_1^2+(\lambda_1^2+\lambda_3^2)x_2^2+(\lambda_1^2+\lambda_2^2)x_3^2+(\lambda_1^2+\lambda_4^2)y_1^2+(\lambda_2^2+\lambda_4^2)y_2^2+(\lambda_3^2+\lambda_4^2)y_3^2\\
&+(\lambda_1^2+\lambda_5^2)z_1^2+(\lambda_2^2+\lambda_5^2)z_2^2+(\lambda_3^2+\lambda_5^2)z_3^2+(\lambda_4^2+\lambda_5^2)z_4^2]
\end{align*}
and respectively
\begin{align*}
K_2(M)&=-\frac{1}{4}\operatorname{Trace}\sum\limits_{p=0}^5 J^{p}MJ^{5-p}\Omega\\
&=\frac{1}{2}[(\lambda_2^4+\lambda_3^4+\lambda_2^2\lambda_3^2)x_1^2+(\lambda_1^4+\lambda_3^4+\lambda_1^2\lambda_3^2)x_2^2+(\lambda_1^4+\lambda_2^4+\lambda_1^2\lambda_2^2)x_3^2+\\
&+(\lambda_1^4+\lambda_4^4+\lambda_1^2\lambda_4^2)y_1^2+(\lambda_2^4+\lambda_4^4+\lambda_2^2\lambda_4^2)y_2^2+(\lambda_3^4+\lambda_4^4+\lambda_3^2\lambda_4^2)y_3^2+\\
&+(\lambda_1^4+\lambda_5^4+\lambda_1^2\lambda_5^2)z_1^2+(\lambda_2^4+\lambda_5^4+\lambda_2^2\lambda_5^2)z_2^2+(\lambda_3^4+\lambda_5^4+\lambda_3^2\lambda_5^2)z_3^2+(\lambda_4^4+\lambda_5^4+\lambda_4^2\lambda_5^2)z_4^2],
\end{align*}
which are, more precisely, the Mishchenko's integrals of order 4 and respectively 6.\\
Thus Mishchenko's method provides two constants of motion (the above constants $K_1$ and $K_2$), which add to the Hamiltonian $H$. As mentioned in \cite{Ratiu80}, the number of independent constants of motion generated by this method is half of the dimension of the adjoint orbit only for $\mathfrak{so}(4)$. It follows that using Mishchenko's method we cannot obtain other independent constants of motion for the rigid body on $\mathfrak{so}(5)$ except the already mentioned ones.\\
In what follows, we will analyze the results of Manakov's method for finding constants of motion (\cite{Manakov76}, \cite{Ratiu80}, \cite{MoPi}), which are to be found as coefficients of the powers of $\gamma$ in the expansion of 
$$\frac{1}{2r}\operatorname{Trace}(M+\gamma J^2)^r$$
for $r=2,3,4,5$.\\
For $r=2$ we obtain the first Casimir, $C_1$. For $r=3$ we obtain the constant of motion $K_1$. For $r=4$ we obtain the constant of motion $K_2$ and the second Casimir, $C_2$. Finally, for $r=5$ we obtain a constant of motion which can also be obtained by Mishchenko's method for order 8, but is dependent of the previously obtained ones, and the new functionally independent constant of motion 
$$K_3(M)=\frac{1}{10}\left(\sum\limits_{j=1}^5T_j\lambda_j^2\right),$$
where
\begin{align*}
T_1&=(x_3z_2-x_2z_3-y_1z_4)^2+(x_2y_3-x_3y_2-z_1z_4)^2+(x_1x_3-y_1y_3-z_1z_3)^2+\\
&+(x_1x_2-y_1y_2-z_1z_2)^2+(x_2^2+x_3^2+y_1^2+z_1^2)^2;
\end{align*} 
\begin{align*}
T_2&=(x_1z_3-x_3z_1-y_2z_4)^2+(x_3y_1-x_1y_3-z_2z_4)^2+(x_2x_3-y_2y_3-z_2z_3)^2+\\
&+(x_1x_2-y_1y_2-z_1z_2)^2+(x_1^2+x_3^2+y_2^2+z_2^2)^2;
\end{align*}
\begin{align*}
T_3&=(x_2z_1-x_1z_2-y_3z_4)^2+(x_1y_2-x_2y_1-z_3z_4)^2+(x_1x_3-y_1y_3-z_1z_3)^2+\\
&+(x_2x_3-y_2y_3-z_2z_3)^2+(x_1^2+x_2^2+y_3^2+z_3^2)^2;
\end{align*}
\begin{align*}
T_4&=(y_1z_1+y_2z_2+y_3z_3)^2+(x_3y_1-x_1y_3-z_2z_4)^2+(x_2y_3-x_3y_2-z_1z_4)^2+\\
&+(x_1y_2-x_2y_1-z_3z_4)^2+(y_1^2+y_2^2+y_3^2+z_4^2)^2;
\end{align*}
\begin{align*}
T_5&=(x_3z_2-x_2z_3-y_1z_4)^2+(x_1z_3-x_3z_1-y_2z_4)^2+(x_2z_1-x_1z_2-y_3z_4)^2+\\
&+(y_1z_1+y_2z_2+y_3z_3)^2+(z_1^2+z_2^2+z_3^2+z_4^2)^2.
\end{align*}
We find, in the particular case of $\mathfrak{so}(5)$, the result from \cite{Uga}, \cite{MoPi}, that the Mishchenko integrals of even order ($m_2=C_1$, $m_4=K_1$, $m_6=K_2$) are a subsystem of the Manakov integrals.

We obtained the following list of functionally independent and Poisson commuting constants of motion for \eqref{soM}: $H,C_1,C_2,K_1,K_2,K_3$, which proves that our system is completely integrable \cite{Manakov76}.

Additionally, as proved in Section 2, the system \eqref{soM} admits the following "generators" integrals of motion:
$$F_1(M)=\frac{x_2^2}{\lambda_1^2-\lambda_3^2}+\frac{x_3^2}{\lambda_1^2-\lambda_2^2}+\frac{y_1^2}{\lambda_1^2-\lambda_4^2}+\frac{z_1^2}{\lambda_1^2-\lambda_5^2};$$
$$F_2(M)=\frac{x_1^2}{\lambda_2^2-\lambda_3^2}+\frac{x_3^2}{\lambda_2^2-\lambda_1^2}+\frac{y_2^2}{\lambda_2^2-\lambda_4^2}+\frac{z_2^2}{\lambda_2^2-\lambda_5^2};$$
$$F_3(M)=\frac{x_1^2}{\lambda_3^2-\lambda_2^2}+\frac{x_2^2}{\lambda_3^2-\lambda_1^2}+\frac{y_3^2}{\lambda_3^2-\lambda_4^2}+\frac{z_3^2}{\lambda_3^2-\lambda_5^2};$$
$$F_4(M)=\frac{y_1^2}{\lambda_4^2-\lambda_1^2}+\frac{y_2^2}{\lambda_4^2-\lambda_2^2}+\frac{y_3^2}{\lambda_4^2-\lambda_3^2}+\frac{z_4^2}{\lambda_4^2-\lambda_5^2};$$
$$F_5(M)=\frac{z_1^2}{\lambda_5^2-\lambda_1^2}+\frac{z_2^2}{\lambda_5^2-\lambda_2^2}+\frac{z_3^2}{\lambda_5^2-\lambda_3^2}+\frac{z_4^2}{\lambda_5^2-\lambda_4^2}.$$

We continue by presenting the list of factorizations of the characteristic polynomials for the linearized equations on the tangent space to the orbit, corresponding to the equilibria from 
$\mathfrak{t}_k\cap \operatorname{Orb}_{c_1;c_2}.$\\
For the equilibrium $M_{a,b}^1$ the eigenvalues are the roots of equations of the following form:
$$\left.\begin{array}{l}
U_1t^2+U_1'=0;~~V_1t^2+V_1'=0;~~W_1t^4+W_1't^2+W_1''=0;\\
U_1=(\lambda_3+\lambda_4)^2(\lambda_3+\lambda_5)(\lambda_4+\lambda_5);~~U_1'=b^2(\lambda_3-\lambda_5)(\lambda_4-\lambda_5);\\
V_1=(\lambda_1+\lambda_2)^2(\lambda_1+\lambda_5)(\lambda_2+\lambda_5);~~V_1'=a^2(\lambda_1-\lambda_5)(\lambda_2-\lambda_5);\\
W_1=(\lambda_1+\lambda_2)^4(\lambda_3+\lambda_4)^4(\lambda_1+\lambda_3)(\lambda_1+\lambda_4)(\lambda_2+\lambda_3)(\lambda_2+\lambda_4);\\
W_1''=(\lambda_1-\lambda_3)(\lambda_1-\lambda_4)(\lambda_2-\lambda_3)(\lambda_2-\lambda_4)[a^2(\lambda_3+\lambda_4)^2-b^2(\lambda_1+\lambda_2)^2]^2.
\end{array}\right.$$
For $M_{a,b}^2$ the eigenvalues are the roots of equations of the following form:
$$\left.\begin{array}{l}
U_2t^2+U_2'=0;~~V_2t^2+V_2'=0;~~W_2t^4+W_2't^2+W_2''=0;\\
U_2=(\lambda_2+\lambda_5)^2(\lambda_1+\lambda_2)(\lambda_1+\lambda_5);~~U_2'=b^2(\lambda_1-\lambda_2)(\lambda_1-\lambda_5);\\
V_2=(\lambda_3+\lambda_4)^2(\lambda_1+\lambda_3)(\lambda_1+\lambda_4);~~V_2'=a^2(\lambda_1-\lambda_3)(\lambda_1-\lambda_4);\\
W_2=(\lambda_2+\lambda_5)^4(\lambda_3+\lambda_4)^4(\lambda_2+\lambda_3)(\lambda_2+\lambda_4)(\lambda_3+\lambda_5)(\lambda_4+\lambda_5);\\
W_2''=(\lambda_2-\lambda_3)(\lambda_2-\lambda_4)(\lambda_3-\lambda_5)(\lambda_4-\lambda_5)[a^2(\lambda_2+\lambda_5)^2-b^2(\lambda_3+\lambda_4)^2]^2.
\end{array}\right.$$
For $M_{a,b}^3$ the eigenvalues are the roots of equations of the following form:
$$\left.\begin{array}{l}
U_3t^2+U_3'=0;~~V_3t^2+V_3'=0;~~W_3t^4+W_3't^2+W_3''=0;\\
U_3=(\lambda_1+\lambda_5)^2(\lambda_1+\lambda_2)(\lambda_2+\lambda_5);~~U_3'=-b^2(\lambda_1-\lambda_2)(\lambda_2-\lambda_5);\\
V_3=(\lambda_3+\lambda_4)^2(\lambda_2+\lambda_3)(\lambda_2+\lambda_4);~~V_3'=a^2(\lambda_2-\lambda_3)(\lambda_2-\lambda_4);\\
W_3=(\lambda_1+\lambda_5)^4(\lambda_3+\lambda_4)^4(\lambda_1+\lambda_3)(\lambda_1+\lambda_4)(\lambda_3+\lambda_5)(\lambda_4+\lambda_5);\\
W_3''=(\lambda_1-\lambda_3)(\lambda_1-\lambda_4)(\lambda_3-\lambda_5)(\lambda_4-\lambda_5)[a^2(\lambda_1+\lambda_5)^2-b^2(\lambda_3+\lambda_4)^2]^2.
\end{array}\right.$$
For $M_{a,b}^4$ the eigenvalues are the roots of equations of the following form:
$$\left.\begin{array}{l}
U_4t^2+U_4'=0;~~V_4t^2+V_4'=0;~~W_4t^4+W_4't^2+W_4''=0;\\
U_4=(\lambda_1+\lambda_5)^2(\lambda_1+\lambda_3)(\lambda_3+\lambda_5);~~U_4'=-b^2(\lambda_1-\lambda_3)(\lambda_3-\lambda_5);\\
V_4=(\lambda_2+\lambda_4)^2(\lambda_2+\lambda_3)(\lambda_3+\lambda_4);~~V_4'=-a^2(\lambda_2-\lambda_3)(\lambda_3-\lambda_4);\\
W_4=(\lambda_1+\lambda_5)^4(\lambda_2+\lambda_4)^4(\lambda_1+\lambda_2)(\lambda_1+\lambda_4)(\lambda_2+\lambda_5)(\lambda_4+\lambda_5);\\
W_4''=(\lambda_1-\lambda_2)(\lambda_1-\lambda_4)(\lambda_2-\lambda_5)(\lambda_4-\lambda_5)[a^2(\lambda_1+\lambda_5)^2-b^2(\lambda_2+\lambda_4)^2]^2.
\end{array}\right.$$
For $M_{a,b}^5$ the eigenvalues are the roots of equations of the following form:
$$\left.\begin{array}{l}
U_5t^2+U_5'=0;~~V_5t^2+V_5'=0;~~W_5t^4+W_5't^2+W_5''=0;\\
U_5=(\lambda_1+\lambda_5)^2(\lambda_1+\lambda_4)(\lambda_4+\lambda_5);~~U_5'=-b^2(\lambda_1-\lambda_4)(\lambda_4-\lambda_5);\\
V_5=(\lambda_2+\lambda_3)^2(\lambda_2+\lambda_4)(\lambda_3+\lambda_4);~~V_5'=a^2(\lambda_2-\lambda_4)(\lambda_3-\lambda_4);\\
W_5=(\lambda_1+\lambda_5)^4(\lambda_2+\lambda_3)^4(\lambda_1+\lambda_2)(\lambda_1+\lambda_3)(\lambda_2+\lambda_5)(\lambda_3+\lambda_5);\\
W_5''=(\lambda_1-\lambda_2)(\lambda_1-\lambda_3)(\lambda_2-\lambda_5)(\lambda_3-\lambda_5)[a^2(\lambda_1+\lambda_5)^2-b^2(\lambda_2+\lambda_3)^2]^2.
\end{array}\right.$$
For $M_{a,b}^6$ the eigenvalues are the roots of equations of the following form:
$$\left.\begin{array}{l}
U_6t^2+U_1'=0;~~V_6t^2+V_6'=0;~~W_6t^4+W_6't^2+W_6''=0;\\
U_6=(\lambda_2+\lambda_4)^2(\lambda_2+\lambda_5)(\lambda_4+\lambda_5);~~U_6'=b^2(\lambda_2-\lambda_5)(\lambda_4-\lambda_5);\\
V_6=(\lambda_1+\lambda_3)^2(\lambda_1+\lambda_5)(\lambda_3+\lambda_5);~~V_6'=a^2(\lambda_1-\lambda_5)(\lambda_3-\lambda_5);\\
W_6=(\lambda_1+\lambda_3)^4(\lambda_2+\lambda_4)^4(\lambda_1+\lambda_2)(\lambda_1+\lambda_4)(\lambda_2+\lambda_3)(\lambda_3+\lambda_4);\\
W_6''=-(\lambda_1-\lambda_2)(\lambda_1-\lambda_4)(\lambda_2-\lambda_3)(\lambda_3-\lambda_4)[a^2(\lambda_2+\lambda_4)^2-b^2(\lambda_1+\lambda_3)^2]^2.
\end{array}\right.$$
For $M_{a,b}^7$ the eigenvalues are the roots of equations of the following form:
$$\left.\begin{array}{l}
U_7t^2+U_7'=0;~~V_7t^2+V_7'=0;~~W_7t^4+W_7't^2+W_7''=0;\\
U_7=(\lambda_3+\lambda_5)^2(\lambda_1+\lambda_3)(\lambda_1+\lambda_5);~~U_7'=b^2(\lambda_1-\lambda_3)(\lambda_1-\lambda_5);\\
V_7=(\lambda_2+\lambda_4)^2(\lambda_1+\lambda_4)(\lambda_1+\lambda_2);~~V_7'=a^2(\lambda_1-\lambda_2)(\lambda_1-\lambda_4);\\
W_7=(\lambda_2+\lambda_4)^4(\lambda_3+\lambda_5)^4(\lambda_2+\lambda_3)(\lambda_2+\lambda_5)(\lambda_3+\lambda_4)(\lambda_4+\lambda_5);\\
W_7''=-(\lambda_2-\lambda_3)(\lambda_2-\lambda_5)(\lambda_3-\lambda_4)(\lambda_4-\lambda_5)[a^2(\lambda_3+\lambda_5)^2-b^2(\lambda_2+\lambda_4)^2]^2.
\end{array}\right.$$
For $M_{a,b}^8$ the eigenvalues are the roots of equations of the following form:
$$\left.\begin{array}{l}
U_8t^2+U_1'=0;~~V_8t^2+V_8'=0;~~W_8t^4+W_8't^2+W_8''=0;\\
U_8=(\lambda_4+\lambda_5)^2(\lambda_1+\lambda_4)(\lambda_1+\lambda_5);~~U_8'=b^2(\lambda_1-\lambda_4)(\lambda_1-\lambda_5);\\
V_8=(\lambda_2+\lambda_3)^2(\lambda_1+\lambda_2)(\lambda_1+\lambda_3);~~V_8'=a^2(\lambda_1-\lambda_2)(\lambda_1-\lambda_3);\\
W_8=(\lambda_2+\lambda_3)^4(\lambda_4+\lambda_5)^4(\lambda_2+\lambda_4)(\lambda_2+\lambda_5)(\lambda_3+\lambda_4)(\lambda_3+\lambda_5);\\
W_8''=(\lambda_2-\lambda_4)(\lambda_2-\lambda_5)(\lambda_3-\lambda_4)(\lambda_3-\lambda_5)[a^2(\lambda_4+\lambda_5)^2-b^2(\lambda_2+\lambda_3)^2]^2.
\end{array}\right.$$
For $M_{a,b}^9$ the eigenvalues are the roots of equations of the following form:
$$\left.\begin{array}{l}
U_9t^2+U_1'=0;~~V_9t^2+V_9'=0;~~W_9t^4+W_9't^2+W_9''=0;\\
U_9=(\lambda_1+\lambda_4)^2(\lambda_1+\lambda_5)(\lambda_4+\lambda_5);~~U_9'=b^2(\lambda_1-\lambda_5)(\lambda_4-\lambda_5);\\
V_9=(\lambda_2+\lambda_3)^2(\lambda_2+\lambda_5)(\lambda_3+\lambda_5);~~V_9'=a^2(\lambda_2-\lambda_5)(\lambda_3-\lambda_5);\\
W_9=(\lambda_1+\lambda_4)^4(\lambda_2+\lambda_3)^4(\lambda_1+\lambda_2)(\lambda_1+\lambda_3)(\lambda_2+\lambda_4)(\lambda_3+\lambda_4);\\
W_9''=(\lambda_1-\lambda_2)(\lambda_1-\lambda_3)(\lambda_2-\lambda_4)(\lambda_3-\lambda_4)[a^2(\lambda_1+\lambda_4)^2-b^2(\lambda_2+\lambda_3)^2]^2.
\end{array}\right.$$
For $M_{a,b}^{10}$ the eigenvalues are the roots of equations of the following form:
$$\left.\begin{array}{l}
U_{10}t^2+U_1'=0;~~V_{10}t^2+V_{10}'=0;~~W_{10}t^4+W_{10}'t^2+W_{10}''=0;\\
U_{10}=(\lambda_2+\lambda_5)^2(\lambda_2+\lambda_4)(\lambda_4+\lambda_5);~~U_{10}'=-b^2(\lambda_2-\lambda_4)(\lambda_4-\lambda_5);\\
V_{10}=(\lambda_1+\lambda_3)^2(\lambda_1+\lambda_4)(\lambda_3+\lambda_4);~~V_{10}'=a^2(\lambda_1-\lambda_4)(\lambda_3-\lambda_4);\\
W_{10}=(\lambda_1+\lambda_3)^4(\lambda_2+\lambda_5)^4(\lambda_1+\lambda_2)(\lambda_1+\lambda_5)(\lambda_2+\lambda_3)(\lambda_3+\lambda_5);\\
W_{10}''=-(\lambda_1-\lambda_2)(\lambda_1-\lambda_5)(\lambda_2-\lambda_3)(\lambda_3-\lambda_5)[a^2(\lambda_2+\lambda_5)^2-b^2(\lambda_1+\lambda_3)^2]^2.
\end{array}\right.$$
For $M_{a,b}^{11}$ the eigenvalues are the roots of equations of the following form:
$$\left.\begin{array}{l}
U_{11}t^2+U_1'=0;~~V_{11}t^2+V_{11}'=0;~~W_{11}t^4+W_{11}'t^2+W_{11}''=0;\\
U_{11}=(\lambda_2+\lambda_5)^2(\lambda_2+\lambda_3)(\lambda_3+\lambda_5);~~U_{11}'=-b^2(\lambda_2-\lambda_3)(\lambda_3-\lambda_5);\\
V_{11}=(\lambda_1+\lambda_4)^2(\lambda_1+\lambda_3)(\lambda_3+\lambda_4);~~V_{11}'=-a^2(\lambda_1-\lambda_3)(\lambda_3-\lambda_4);\\
W_{11}=(\lambda_1+\lambda_4)^4(\lambda_2+\lambda_5)^4(\lambda_1+\lambda_2)(\lambda_1+\lambda_5)(\lambda_2+\lambda_4)(\lambda_4+\lambda_5);\\
W_{11}''=-(\lambda_1-\lambda_2)(\lambda_1-\lambda_5)(\lambda_2-\lambda_4)(\lambda_4-\lambda_5)[a^2(\lambda_2+\lambda_5)^2-b^2(\lambda_1+\lambda_4)^2]^2.
\end{array}\right.$$
For $M_{a,b}^{12}$ the eigenvalues are the roots of equations of the following form:
$$\left.\begin{array}{l}
U_{12}t^2+U_1'=0;~~V_{12}t^2+V_{12}'=0;~~W_{12}t^4+W_{12}'t^2+W_{12}''=0;\\
U_{12}=(\lambda_4+\lambda_5)^2(\lambda_3+\lambda_4)(\lambda_3+\lambda_5);~~U_{12}'=b^2(\lambda_3-\lambda_4)(\lambda_3-\lambda_5);\\
V_{12}=(\lambda_1+\lambda_2)^2(\lambda_1+\lambda_3)(\lambda_2+\lambda_3);~~V_{12}'=a^2(\lambda_1-\lambda_3)(\lambda_2-\lambda_3);\\
W_{12}=(\lambda_1+\lambda_2)^4(\lambda_4+\lambda_5)^4(\lambda_1+\lambda_4)(\lambda_1+\lambda_5)(\lambda_2+\lambda_4)(\lambda_2+\lambda_5);\\
W_{12}''=(\lambda_1-\lambda_4)(\lambda_1-\lambda_5)(\lambda_2-\lambda_4)(\lambda_2-\lambda_5)[a^2(\lambda_4+\lambda_5)^2-b^2(\lambda_1+\lambda_2)^2]^2.
\end{array}\right.$$
For $M_{a,b}^{13}$ the eigenvalues are the roots of equations of the following form:
$$\left.\begin{array}{l}
U_{13}t^2+U_1'=0;~~V_{13}t^2+V_{13}'=0;~~W_{13}t^4+W_{13}'t^2+W_{13}''=0;\\
U_{13}=(\lambda_3+\lambda_5)^2(\lambda_3+\lambda_4)(\lambda_4+\lambda_5);~~U_{13}'=-b^2(\lambda_3-\lambda_4)(\lambda_4-\lambda_5);\\
V_{13}=(\lambda_1+\lambda_2)^2(\lambda_1+\lambda_4)(\lambda_2+\lambda_4);~~V_{13}'=a^2(\lambda_1-\lambda_4)(\lambda_2-\lambda_4);\\
W_{13}=(\lambda_1+\lambda_2)^4(\lambda_3+\lambda_5)^4(\lambda_1+\lambda_3)(\lambda_1+\lambda_5)(\lambda_2+\lambda_3)(\lambda_2+\lambda_5);\\
W_{13}''=(\lambda_1-\lambda_3)(\lambda_1-\lambda_5)(\lambda_2-\lambda_3)(\lambda_2-\lambda_5)[a^2(\lambda_3+\lambda_5)^2-b^2(\lambda_1+\lambda_2)^2]^2.
\end{array}\right.$$
For $M_{a,b}^{14}$ the eigenvalues are the roots of equations of the following form:
$$\left.\begin{array}{l}
U_{14}t^2+U_1'=0;~~V_{14}t^2+V_{14}'=0;~~W_{14}t^4+W_{14}'t^2+W_{14}''=0;\\
U_{14}=(\lambda_3+\lambda_5)^2(\lambda_2+\lambda_3)(\lambda_2+\lambda_5);~~U_{14}'=b^2(\lambda_2-\lambda_3)(\lambda_2-\lambda_5);\\
V_{14}=(\lambda_1+\lambda_4)^2(\lambda_1+\lambda_2)(\lambda_2+\lambda_4);~~V_{14}'=-a^2(\lambda_1-\lambda_2)(\lambda_2-\lambda_4);\\
W_{14}=(\lambda_1+\lambda_4)^4(\lambda_3+\lambda_5)^4(\lambda_1+\lambda_3)(\lambda_1+\lambda_5)(\lambda_3+\lambda_4)(\lambda_4+\lambda_5);\\
W_{14}''=-(\lambda_1-\lambda_3)(\lambda_1-\lambda_5)(\lambda_3-\lambda_4)(\lambda_4-\lambda_5)[a^2(\lambda_3+\lambda_5)^2-b^2(\lambda_1+\lambda_4)^2]^2.
\end{array}\right.$$
For $M_{a,b}^{15}$ the eigenvalues are the roots of equations of the following form:
$$\left.\begin{array}{l}
U_{15}t^2+U_1'=0;~~V_{15}t^2+V_{15}'=0;~~W_{15}t^4+W_{15}'t^2+W_{15}''=0;\\
U_{15}=(\lambda_4+\lambda_5)^2(\lambda_2+\lambda_5)(\lambda_2+\lambda_4);~~U_{15}'=b^2(\lambda_2-\lambda_4)(\lambda_2-\lambda_5);\\
V_{15}=(\lambda_1+\lambda_3)^2(\lambda_1+\lambda_2)(\lambda_2+\lambda_3);~~V_{15}'=-a^2(\lambda_1-\lambda_2)(\lambda_2-\lambda_3);\\
W_{15}=(\lambda_1+\lambda_3)^4(\lambda_4+\lambda_5)^4(\lambda_1+\lambda_4)(\lambda_1+\lambda_5)(\lambda_3+\lambda_4)(\lambda_3+\lambda_5);\\
W_{15}''=(\lambda_1-\lambda_4)(\lambda_1-\lambda_5)(\lambda_3-\lambda_4)(\lambda_3-\lambda_5)[a^2(\lambda_4+\lambda_5)^2-b^2(\lambda_1+\lambda_3)^2]^2.
\end{array}\right.$$
Without loss of generality, we can choose an ordering for $\lambda_i$'s, namely
$$\lambda_1>\lambda_2>\lambda_3>\lambda_4>\lambda_5.$$

\medskip

Since for the equilibria $\mathfrak{t}_3 \cap \operatorname{Orb}_{c_1;c_2}$, $\mathfrak{t}_4 \cap \operatorname{Orb}_{c_1;c_2}$, $\mathfrak{t}_5 \cap \operatorname{Orb}_{c_1;c_2}$, $\mathfrak{t}_{10} \cap \operatorname{Orb}_{c_1;c_2}$, $\mathfrak{t}_{11} \cap \operatorname{Orb}_{c_1;c_2}$, $\mathfrak{t}_{13} \cap \operatorname{Orb}_{c_1;c_2}$, $\mathfrak{t}_{14} \cap \operatorname{Orb}_{c_1;c_2}$ and $\mathfrak{t}_{15} \cap \operatorname{Orb}_{c_1;c_2}$ at least one of the eigenvalues is real and strictly positive we have nonlinear {\it instability} for these equilibria. 

\medskip

It is immediate to prove that if
$$a^2(\lambda_2+\lambda_4)^2\not= b^2(\lambda_1+\lambda_3)^2$$
is equivalent from Eqs. \ref{values_a_b} with
\begin{equation}
\label{special1}
c_1^2\left[(\lambda_2+\lambda_4)^4+(\lambda_1+\lambda_3)^4\right]\not= c_2\left[(\lambda_2+\lambda_4)^2+(\lambda_1+\lambda_3)^2\right]^2,
\end{equation}
and in this case we have $W_6''<0$. 
Since $W_6>0$ and $W_6''<0$, it follows that if the condition \eqref{special1} holds, then the equilibrium $M_{a,b}^6$ from $\mathfrak{t}_6 \cap \operatorname{Orb}_{c_1;c_2}$ (together with the equilibria $M_{-a,b}^6$, $M_{a,-b}^6$, $M_{-a,-b}^6$) are also nonlinear {\it unstable}.
 
On the other hand, since
$$b^2(\lambda_2+\lambda_4)^2< a^2(\lambda_1+\lambda_3)^2,$$
(because $a>b$ and $\lambda_1+\lambda_3>\lambda_2+\lambda_4$)
it follows that the equilibrium $M_{b,a}^6$ from $\mathfrak{t}_6 \cap \operatorname{Orb}_{c_1;c_2}$ (together with the equilibria $M_{b,-a}^6$, $M_{-b,a}^6$, $M_{-b,-a}^6$) are nonlinear {\it unstable}, regardless if condition \eqref{special1} holds or not.

\medskip

A similar discussion holds for the equilibria from $\mathfrak{t}_7\cap \operatorname{Orb}_{c_1;c_2}$.

\medskip

We continue with the study of the remaining families of equilibria, which is $\mathfrak{t}_1\cap \operatorname{Orb}_{c_1;c_2}$, $\mathfrak{t}_2\cap \operatorname{Orb}_{c_1;c_2}$, $\mathfrak{t}_8\cap \operatorname{Orb}_{c_1;c_2}$, $\mathfrak{t}_9\cap \operatorname{Orb}_{c_1;c_2}$ and $\mathfrak{t}_{12}\cap \operatorname{Orb}_{c_1;c_2}$. 

For start, let us notice that energy methods using the "canonical" integrals of motion (the Hamiltonian and/or the constants of motion $K_1,K_2,K_3$) work only for a few of the above equilibria. But choosing particular linear combinations of some "generators" integrals of motion introduced in Section 2 proves to be very effective in solving the stability problem with energy methods.

To begin with, we note that the method involving linearization used above is inconclusive for the stability of the $M_{a,b}^1$ equilibrium. 
Therefore, we will use Arnold's method \cite{Arnold65}, which is equivalent with the other energy methods \cite{BiPu07}. 

Consider the smooth function $G_{mn}\in
C^{\infty}(\mathfrak{so}(5),\mathbb{R})$, where $m,n$ are real numbers,
$$G_{mn}(M)=F_1(M)+F_5(M)+m C_1(M)+n C_2(M).$$
Choosing $m,n$ such that $\mathbf{d}G_{mn}(M_{a,b}^1)=0$, namely
$$m=\frac{2b^2}{(\lambda_1^2-\lambda_2^2)(a^2-b^2)},~~~n=-\frac{2}{(\lambda_1^2-\lambda_2^2)(a^2-b^2)}$$
and taking into account that
$$W:=\ker dC_1(M_{a,b}^1)\cap \ker dC_2(M_{a,b}^1)=\operatorname{span} (E_1,E_2,E_4,E_5,E_7,E_8,E_9,E_{10})$$
we obtain the determinants associated to all upper-left submatrices of the
Hessian $\mathbf{d}^2G_{mn}(M_{a,b}^1)|_{W\times W}$ as follows:
\begin{align*}
D_1&=-\frac{2a^2}{(\lambda_1^2-\lambda_2^2)(a^2-b^2)}<0;\\
D_2&=\frac{4a^2[a^2(\lambda_2^2-\lambda_3^2)+b^2(\lambda_1^2-\lambda_2^2)]}{(\lambda_1^2-\lambda_2^2)^2(\lambda_1^2-\lambda_3^2)(a^2-b^2)^2}>0;\\
D_3&=-\frac{8a^2(\lambda_2^2-\lambda_4^2)[a^2(\lambda_2^2-\lambda_3^2)+b^2(\lambda_1^2-\lambda_2^2)]}{(\lambda_1^2-\lambda_2^2)^3(\lambda_1^2-\lambda_3^2)(\lambda_1^2-\lambda_4^2)(a^2-b^2)^2}<0;\\
D_4&=\frac{16a^4(\lambda_2^2-\lambda_3^2)(\lambda_2^2-\lambda_4^2)}{(\lambda_1^2-\lambda_2^2)^4(\lambda_1^2-\lambda_3^2)(\lambda_1^2-\lambda_4^2)(a^2-b^2)^2}>0;\\
D_5&=-\frac{32a^4(\lambda_2^2-\lambda_4^2)(\lambda_2^2-\lambda_3^2)}{(\lambda_1^2-\lambda_2^2)^5(\lambda_1^2-\lambda_3^2)(\lambda_1^2-\lambda_4^2)(a^2-b^2)^2}<0;\\
D_6&=\frac{64a^4(\lambda_1^2-\lambda_5^2)(\lambda_2^2-\lambda_3^2)(\lambda_2^2-\lambda_4^2)}{(\lambda_1^2-\lambda_2^2)^6(\lambda_1^2-\lambda_3^2)(\lambda_1^2-\lambda_4^2)(\lambda_2^2-\lambda_5^2)(a^2-b^2)^2}>0;\\
D_7&=-\frac{128a^4(\lambda_1^2-\lambda_5^2)(\lambda_2^2-\lambda_3^2)(\lambda_2^2-\lambda_4^2)}{(\lambda_1^2-\lambda_2^2)^6(\lambda_1^2-\lambda_3^2)(\lambda_1^2-\lambda_4^2)(\lambda_2^2-\lambda_5^2)(\lambda_3^2-\lambda_5^2)(a^2-b^2)^2}<0;\\
D_8&=\frac{256a^4(\lambda_1^2-\lambda_5^2)(\lambda_2^2-\lambda_3^2)(\lambda_2^2-\lambda_4^2)}{(\lambda_1^2-\lambda_2^2)^6(\lambda_1^2-\lambda_3^2)(\lambda_1^2-\lambda_4^2)(\lambda_2^2-\lambda_5^2)(\lambda_3^2-\lambda_5^2)(\lambda_4^2-\lambda_5^2)(a^2-b^2)^2}>0.
\end{align*}
Consequently, $\mathbf{d}^2G_{mn}(M_{a,b}^1)|_{W\times W}$ is negative definite, which implies nonlinear {\it stability} for the equilibrium $M_{a,b}^1$. The same computations lead to nonlinear stability for the equilibria $M_{a,-b}^1$, $M_{-a,b}^1$ and $M_{-a,-b}^1$ from $\mathfrak{t}_1\cap \operatorname{Orb}_{c_1;c_2}$.

For the remaining equilibria $M_{b,a}^1$, $M_{-b,a}^1$, $M_{b,-a}^1$ and $M_{-b,-a}^1$ from $\mathfrak{t}_1\cap \operatorname{Orb}_{c_1;c_2}$ a convenient energy function for applying Arnold's method is $G_{mn}(M)=F_4(M)+F_5(M)+m C_1(M)+n C_2(M)$. Thus, all eight equilibria in $\mathfrak{t}_1\cap \operatorname{Orb}_{c_1;c_2}$ are nonlinear {\it stable}.

\medskip

For the equilibria $M_{a,b}^2$, $M_{-a,b}^2$, $M_{a,-b}^2$ and $M_{-a,-b}^2$ from $\mathfrak{t}_2\cap \operatorname{Orb}_{c_1;c_2}$ a convenient energy function for applying Arnold's method is $G_{mn}(M)=F_1(M)-F_4(M)+m C_1(M)+n C_2(M)$. Thus, the equilibria $M_{a,b}^2$, $M_{-a,b}^2$, $M_{a,-b}^2$ and $M_{-a,-b}^2$ are nonlinear {\it stable}.

\medskip

For the equilibria $M_{a,b}^8$, $M_{-a,b}^8$, $M_{a,-b}^8$ and $M_{-a,-b}^8$ from $\mathfrak{t}_8\cap \operatorname{Orb}_{c_1;c_2}$ a convenient energy function for applying Arnold's method is $G_{mn}(M)=F_1(M)+F_5(M)+m C_1(M)+n C_2(M)$. For the equilibria $M_{b,a}^8$, $M_{-b,a}^8$, $M_{b,-a}^8$ and $M_{-b,-a}^8$ from $\mathfrak{t}_8\cap \operatorname{Orb}_{c_1;c_2}$ a convenient energy function for applying Arnold's method is $G_{mn}(M)=F_1(M)+F_2(M)+m C_1(M)+n C_2(M)$. Thus, all eight equilibria in $\mathfrak{t}_8\cap \operatorname{Orb}_{c_1;c_2}$ are nonlinear {\it stable}.

\medskip

For the equilibria $M_{b,a}^9$, $M_{-b,a}^9$, $M_{b,-a}^9$ and $M_{-b,-a}^9$ from $\mathfrak{t}_9\cap \operatorname{Orb}_{c_1;c_2}$ a convenient energy function for applying Arnold's method is $G_{mn}(M)=F_4(M)-F_5(M)+m C_1(M)+n C_2(M)$. Thus, the equilibria $M_{b,a}^9$, $M_{-b,a}^9$, $M_{b,-a}^9$ and $M_{-b,-a}^9$ from $\mathfrak{t}_9\cap \operatorname{Orb}_{c_1;c_2}$ are nonlinear {\it stable}.

\medskip

For the equilibria $M_{a,b}^{12}$, $M_{-a,b}^{12}$, $M_{a,-b}^{12}$ and $M_{-a,-b}^{12}$ from $\mathfrak{t}_{12}\cap \operatorname{Orb}_{c_1;c_2}$ a convenient energy function for applying Arnold's method is $G_{mn}(M)=F_1(M)-F_3(M)+m C_1(M)+n C_2(M)$. For the equilibria $M_{b,a}^{12}$, $M_{-b,a}^{12}$, $M_{b,-a}^{12}$ and $M_{-b,-a}^{12}$ from $\mathfrak{t}_{12}\cap \operatorname{Orb}_{c_1;c_2}$ a convenient energy function for applying Arnold's method is $G_{mn}(M)=F_3(M)-F_5(M)+m C_1(M)+n C_2(M)$. Thus, all eight equilibria in $\mathfrak{t}_{12}\cap \operatorname{Orb}_{c_1;c_2}$ are nonlinear {\it stable}.

\medskip

For the equilibria $M_{b,a}^2$, $M_{-b,a}^2$, $M_{b,-a}^2$ and $M_{-b,-a}^2$ from $\mathfrak{t}_2\cap \operatorname{Orb}_{c_1;c_2}$ and for the equilibria $M_{a,b}^9$, $M_{-a,b}^9$, $M_{a,-b}^9$ and $M_{-a,-b}^9$ from $\mathfrak{t}_9\cap \operatorname{Orb}_{c_1;c_2}$ the stability problem remains open and it is likely that a bifurcation phenomenon occurs. Such a phenomenon appears for some equilibria in the case of 
the $\mathfrak{so}(4)$ free rigid body and is extensively studied in \cite{noi}. 

\medskip

We proved the following result:
\begin{thm}
\begin{itemize}
\item[{\rm (i)}] The equilibria from $\mathfrak{t}_3 \cap \operatorname{Orb}_{c_1;c_2}$, $\mathfrak{t}_4 \cap \operatorname{Orb}_{c_1;c_2}$, $\mathfrak{t}_5 \cap \operatorname{Orb}_{c_1;c_2}$, $\mathfrak{t}_{10} \cap \operatorname{Orb}_{c_1;c_2}$, $\mathfrak{t}_{11} \cap \operatorname{Orb}_{c_1;c_2}$, $\mathfrak{t}_{13} \cap \operatorname{Orb}_{c_1;c_2}$, $\mathfrak{t}_{14} \cap \operatorname{Orb}_{c_1;c_2}$ and $\mathfrak{t}_{15} \cap \operatorname{Orb}_{c_1;c_2}$ are unstable.
\item[{\rm (ii)}] The equilibria $M_{b,a}^6$, $M_{b,-a}^6$, $M_{-b,a}^6$ and $M_{-b,-a}^6$ from $\mathfrak{t}_6 \cap \operatorname{Orb}_{c_1;c_2}$ and respectively the equilibria $M_{b,a}^7$, $M_{b,-a}^7$, $M_{-b,a}^7$ and $M_{-b,-a}^7$ from $\mathfrak{t}_7 \cap \operatorname{Orb}_{c_1;c_2}$ are unstable.
\item[{\rm (iii)}] a) If the condition
$$c_1^2\left[(\lambda_1+\lambda_3)^4+(\lambda_2+\lambda_4)^4\right]\not= c_2\left[(\lambda_1+\lambda_3)^2+(\lambda_2+\lambda_4)^2\right]^2$$
holds, then the equilibria $M_{a,b}^6$, $M_{-a,b}^6$, $M_{a,-b}^6$ and $M_{-a,-b}^6$ from $\mathfrak{t}_6 \cap \operatorname{Orb}_{c_1;c_2}$ are unstable.\\
b) If the condition
$$c_1^2\left[(\lambda_2+\lambda_4)^4+(\lambda_3+\lambda_5)^4\right]\not= c_2\left[(\lambda_2+\lambda_4)^2+(\lambda_3+\lambda_5)^2\right]^2$$
holds, then the equilibria $M_{a,b}^7$, $M_{-a,b}^7$, $M_{a,-b}^7$ and $M_{-a,-b}^7$ from $\mathfrak{t}_7 \cap \operatorname{Orb}_{c_1;c_2}$ are unstable.
\item[{\rm (iv)}] The equilibria from $\mathfrak{t}_1 \cap \operatorname{Orb}_{c_1;c_2}$, $\mathfrak{t}_8 \cap \operatorname{Orb}_{c_1;c_2}$ and $\mathfrak{t}_{12} \cap \operatorname{Orb}_{c_1;c_2}$ are nonlinear stable.
\item[{\rm (v)}] The equilibria $M_{a,b}^2$, $M_{-a,b}^2$, $M_{a,-b}^2$ and $M_{-a,-b}^2$ from $\mathfrak{t}_2 \cap \operatorname{Orb}_{c_1;c_2}$ and respectively the equilibria $M_{b,a}^9$, $M_{-b,a}^9$, $M_{b,-a}^9$ and $M_{-b,-a}^9$ from $\mathfrak{t}_9\cap \operatorname{Orb}_{c_1;c_2}$ are nonlinear {\it stable}.
\end{itemize}
\end{thm}

\medskip

\noindent {\bf Acknowledgements.} The author acknowledges the very helpful discussions with Professors Tudor Ra\c tiu and Petre Birtea.

\end{document}